\documentstyle [12pt] {article}
\input{amssym.def}
\input{amssym.tex}

\newcommand{\Q}{{\Bbb Q}}
\newcommand{\R}{{\Bbb R}}
\newcommand{\Z}{{\Bbb Z}}

\newcommand{\wreath}{\mbox{ \rm wr }}

\newtheorem{*example}[defi]{*Example}

\addtocontents{ences.txt}{References}

\begin{document}

\title{\bf Unsolved Problems in Ordered and Orderable Groups}
\author{\small V.\,V.~Bludov, A.\,M.\,W.~Glass, V.\,M.~Kopytov
and \framebox{N. Ya. Medvedev}}
\maketitle

Dedicated to Paul F. Conrad, in memoriam.

\bigskip

\setcounter{theorem}{0}

\section*{Introduction}

The study of partially ordered and lattice-ordered groups arose from 
three distinct sources.
\medskip

The first source is functional analysis; in particular, the study of 
real-valued functions on a topological space. Besides the addition of 
functions, it also took into consideration the pointwise  
ordering. Classes of continuous, differentiable, piecewise linear, 
piecewise continuous, and piecewise differentiable functions are especially natural classes that are included in the study. 
For an extensive account of this work, see \cite{Lux1}.
\medskip

The second source came from lectures by Philip Hall in Cambridge and 
A. I. Malcev in Novosibirsk. These focused on finding necessary and 
sufficient conditions for a group to be orderable. This leads to a 
study which is, by its nature, primarily group theoretic. In this 
light, lattice-ordered groups can be seen as groups with operators, 
which connect them up with universal algebra in a natural way.
\medskip

The third source came from generalisations of valuation theory. This 
was initiated by Paul F. Conrad (among others) in his Ph.D. thesis 
under Reinhold Baer. Through his work with his research/post-doctoral 
students John Harvey and W. Charles Holland, this lead to the 
study of the structure theory for abelian lattice-ordered groups and 
later to the structure theory in the non-abelian case as well. This 
attempted to prove similar results to those in group theory and 
centred on finding necessary and sufficient conditions for the 
lattices of certain (convex sublattice) subgroups to be especially 
nice. The success of this program can be seen by the sparsity of 
questions that remain in this topic since 1990.
In contrast to the second approach, Conrad's work focused on the 
importance of the lattice order as opposed to the group operation, and 
it was no surprise that many of the talks at the Conrad Memorial Conference in 2009 reflected 
this. Unfortunately, nothing as useful as the Hirsch-Plotkin radical 
has yet surfaced, a challenge to the next generation of researchers.
One recent possibility is the Conrad soul (subgroup) for right-ordered groups introduced by Andr\'{e}s Navas (see \cite{Na} and \cite{NRC}).
\medskip

The ability to marry the second and third approaches was obtained as a 
result of W. Charles Holland's representation theorem for lattice-ordered groups using permutation groups.
This lead to the study of varieties and quasi-varieties of 
lattice-ordered groups, verbal and marginal subgroups, torsion classes 
and similar concepts of interest to universal algebraists and group 
theorists alike.
\medskip 

To help focus research in the subject, V. M. Kopytov and N. Ya.~Med\-ve\-dev compiled 
a list of problems that had appeared scattered 
through the literature. This appeared in the Nikolai Ya. Medvedev 
Memorial Volume \cite{MMV} (in Russian). Here is an English translation of 
that list provided initially by Medvedev himself in conversations with us. We 
have also added recent new problems, and provided comments where possible,
with a list of solutions where known. The list is along the same lines as the famous Kourovka Note Book \cite{K289} which has been so vital in providing research topics in group theory.

\section{Lattice-ordered groups}

Those problems marked with $^*$ have comments in subsequent subsections.
\medskip

We use $\ell$-group as an abbreviation for lattice-ordered group. As 
is standard, a subdirect product of ordered groups is called a {\it 
residually ordered} $\ell$-group. 
Other names for the same concept are {\it o-approximable} and {\it 
representable}. The class of all residually ordered $\ell$-groups 
forms a variety which will be denoted by ${\cal R}=\;${\scriptsize 
R}$O$. 

The free $\ell$-group on $n$ generators ($n\in \Z_+$) is denoted by 
$F_n$.

A group $G$ is called an {\it Engel group} if for all $x,y\in G$ there 
is $n\in \Z_+$ such that $[x,_n y]=1$ (where $[x,_n y]$ is a shorthand for the 
commutator of $x$ followed by $n$ $y'$s). If $n\in\Z_+$ can be chosen independent
of the pair of elements $x,y\in G$, then the group is called $n$-{\it Engel} (or, in Russian, {\it boundary Engel}).

A variety ${\cal V}$ is said to have {\it finite rank} if there is a positive integer $n$ such that ${\cal V}$ can be defined by equations involving only $x_1,\dots,x_n$, and to be {\it finitely based} if it can be defined by a finite set of equations. 
\medskip

\subsection{Open Problems}

\noindent1.1. If an $\ell$-group $G$ is $\ell$-simple, are all its convex normal 
subgroups trivial?
\medskip

\noindent1.2. Find conditions on an $\ell$-homomorphism for it to 
preserve the ar\-chi\-me\-de\-an property.
\medskip

\noindent1.3. Describe all lattice orders on the free soluble 
$\ell$-group on a finite number $n$ ($>1$) of generators.
\medskip

\noindent1.4. If $H$ is a countable $\ell$-group, is there an abelian 
$\ell$-group $G$ that is isomorphic to $H$ as a lattice?
\medskip

\noindent1.5$^*$. Can the lattice order on a free $\ell$-group be extended 
to a weakly abelian total order?
\medskip

\noindent1.6. For any even positive integer $k$, does there exist 
an $\ell$-group admitting exactly $k$ distinct lattice orders?  (cf. Problem 3.24)
\medskip

\noindent1.7$^*$. (E. C. Weinberg \cite{H283} Problem 2) Does every torsion-free abelian 
group have an archimedean lattice order?
\medskip

\noindent1.8. Find (in explicit form) a variety of $\ell$-groups that has no 
finite basis of laws but every proper subvariety has a finite basis of laws.
\medskip

\noindent1.9. If ${\cal V}_1$ and ${\cal V}_2$ have finite rank, does 
${\cal V}_1{\cal V}_2$? 
\medskip

\noindent1.10$^*$. (A. M. W. Glass, W. Ch. Holland, S. H. McCleary \cite{G284}) 
If ${\cal V}_1$ and ${\cal V}_2$ have finite bases of laws, does 
${\cal V}_1{\cal V}_2$? 
\medskip

\noindent1.12$^*$. Is every Engel $\ell$-group residually ordered?  
(cf. Problems 2.8, 2.9, 2.25, 2.29, 3.4, 3.8, 3.9, 3.18)
\medskip

\noindent1.14$^*$. Describe the idempotents in the semigroup of 
quasi-varieties of $\ell$-groups. What is its cardinality?
\medskip

\noindent1.15. If a quasi-variety ${\frak X}$ is defined by a finite 
set of quasi-identities, then is $\ell$-var(${\frak X}$) finitely based?
\medskip

\noindent1.16$^*$. Is $Th(F_n)$ decidable if $n\neq 1$? 
\medskip

\noindent1.17$^*$. If $n\in \Z_+$ with $n\neq 1$, is the universal theory 
of $F_n$ decidable?
\medskip

\noindent1.18$^*$. (A. M. W. Glass \cite{H283}, Problem 41) If $n$ and $m$ are distinct 
positive integers, are $Th(F_n)$ and $Th(F_m)$ distinct?
\medskip

\noindent1.20$^*$. (P. F. Conrad) If an $\ell$-group $G$ is a subdirect product 
of Ohkuma groups, then is addition defined by its lattice? (cf. Problem 1.22)
\medskip

\noindent1.21. (P. F. Conrad, M. R. Darnel \cite{C20}) Is every $\ell$-group with 
unique addition archimedean?
\medskip

\noindent1.22$^*$. (P. F. Conrad, M. R. Darnel \cite{C20}) If $G$ is a subdirect product 
of Ohkuma groups, does $G$ have unique addition? (cf. Problem 1.20)
\medskip

\noindent1.23. (P. F. Conrad, M. R. Darnel \cite{C20}) Is the class of $\ell$-groups 
with unique addition closed under $\ell$-ideals?
\medskip

\noindent1.26$^*$. (A. M. W. Glass \cite{K289}, Problem 12.12) Is the conjugacy problem for 
nilpotent finitely presented $\ell$-groups soluble?
\medskip

\noindent1.27$^*$. (W. B. Powell and C. Tsinakis \cite{H283} Problem 27(a)) Does any 
non-abelian variety of $\ell$-groups have the amalgamation property?
\medskip

\noindent1.28$^*$. Let $(\Omega,\leq)$ be a totally ordered set. What is the 
quasi-variety of $\ell$-groups generated by $Aut(\Omega,\leq)$?
\medskip

\noindent1.29. (S. H. McCleary \cite{K289} Problem 12.52) If $n\in \Z_+$ with $n\neq 1$, is 
$F_n$ Hopfian (i.e., $F_n$ is not $\ell$-isomorphic to a proper quotient of itself)?
\medskip

\noindent1.30. (S. H. McCleary \cite{K289} Problem 12.53) Do free $\ell$-groups on finitely 
many generators have soluble conjugacy problem?
\medskip

\noindent1.31. (S. H. McCleary \cite{K289} Problem 12.54) Does there exist a normal-valued 
$\ell$-group whose lattice of convex $\ell$-subgroups is not lattice-isomorphic 
to that of any abelian $\ell$-group?
\medskip

\noindent1.32$^*$. What is the cardinality of the set of all subvarieties 
of ${\cal R}$ without covers (in ${\cal R}$)?
\medskip

\noindent1.33. (P. F. Conrad, see \cite{H283} Problem 28 written in by J. Martinez) Is the class of all abelian 
$\ell$-groups which support vector lattices a torsion class?
\medskip

\noindent1.34. (D. A. Van Rie \cite{VR94}) If $G=\Z \wreath \Z$ with the usual 
lattice order, does the quasi-variety of $\ell$-groups generated by $G$ 
cover ${\cal A}$ in the lattice of quasi-varieties of $\ell$-groups?
\medskip

\noindent1.35. (M. A. Bardwell \cite{B499}) Is every $\ell$-group 
$\ell$-isomorphic to $Aut(\Omega,\leq)$ for some {\it partially} ordered set 
$(\Omega,\leq)$ ?
\medskip

\noindent1.36$^*$. (W. B. Powell, C. Tsinakis \cite{P293} and \cite{H283} Problem 15) In what varieties 
of $\ell$-groups does the free product preserve the solubility of the word 
problem?
\medskip

\noindent1.37. (W. B. Powell, C. Tsinakis \cite{P293}) In what varieties of $\ell$-groups is the $\ell$-subgroup generated by the factors always their free product?
\medskip

\noindent1.38. (P. F. Conrad \cite{H283} Problem 29) Are any two scalar products on vector 
lattices conjugate by a continuous $\ell$-automorphism?
\medskip

\noindent1.39. Is the $\ell$-subgroup of the free product of $\ell$-groups generated by commutators always free (as an $\ell$-group)?
\medskip

\noindent1.40$^*$. If every $\ell$-subgroup of an $\ell$-group $G$ is 
closed, then is $G$ discrete?
\medskip

\noindent1.41$^*$. Is an $\ell$-group generated by compact elements iff 
each of its $\ell$-subgroups is closed?
\medskip

\noindent1.42. If $G$ is a p.o. group which is the intersection of 
right orders, then can $G$ be order-embedded in an $\ell$-group with 
the embedding preserving any existing sups and infs?
\medskip

\noindent1.44. (V. I. Pshenichnov) If an $\ell$-group is also a topological 
space with continuous group and lattice operations, then such a group is called a
$t\ell$-{\it group}. Do free $t\ell$-groups exist? If yes, describe 
them.
\medskip

\noindent1.45. (V. I. Pshenichnov) Is a connected $t\ell$-group archimedian?
\medskip

\noindent1.46. (V. I. Pshenichnov) Is a connected $t\ell$-group of finite rank 
isomorphic to $\R^n$ as a group?
\medskip

\noindent1.47. Let a $t\ell$-group contain two closed normal 
$t\ell$-subgroups and be $\ell$-isomorphic to their direct sum. Is the group
topologically isomorphic to their direct sum?

\begin{center}
{\bf New problems}
\end{center}

\noindent1.48. Find a variety of $\ell$-groups that does not have 
finite rank but every proper subvariety has finite rank.
\medskip

\noindent1.49. (P. F. Conrad \cite{H283} Problem 9) Can we recognise the closed convex 
$\ell$-subgroups within the lattice of all convex $\ell$-subgroups of a normal-valued 
$\ell$-group $G$?
\medskip

\noindent1.50. Is the elementary theory of the lattices of $\ell$-ideals of 
$\ell$-groups decidable?
\medskip

\noindent1.51. (V. V. Bludov, A. M. W. Glass \cite{BG08}) Can the free product of 
ordered groups with order-isomorphic amalgamated subgroups ever be 
lattice-orderable if it is not orderable?
\medskip

\noindent1.52.  (V. V. Bludov, A. M. W. Glass \cite{BG06}) Is the variety of 
$\ell$-groups generated by all nilpotent $\ell$-groups finitely based?
\medskip

\noindent 1.53. Study the structures $(G,u)$ where $G$ is an $\ell$-group and $u$ is a constant symbol/element of $G$ satisfying the non-first-order sentence
$(\forall g\in G)(\exists n=n(g)\in \Z_+)(u^{-n}\leq g\leq u^n)$.
These are closely related to algebras arising from many-valued logics.
\medskip

\noindent 1.54. (G. Sabbagh) If finitely presented $\ell$-groups satisfy the same first order sentences, must they be isomorphic? 
\medskip

\noindent 1.55$^*$. (D. Mundici) If finitely presented abelian $\ell$-groups satisfy the same first order sentences, must they be isomorphic? 
\medskip

\noindent 1.56$^*$. (D. Mundici) Develop the theory of abelian unital $\ell$-groups.
In particular, characterise and classify the finitely generated projectives in this class (and the finitely presented objects in the equivalent category of MV-algebras).
\medskip

\noindent 1.57$^*$. (G. Sabbagh) What can one say about quasi-finitely axiomatisable (QFA) $\ell$-groups, abelian $\ell$-groups and totally ordered groups. 
(A finitely generated model $G$ of an algebraic theory $T$ is said to be QFA
if there is a sentence $\sigma$ which holds in $G$ such that every finitely generated model of $T$ which satisfies $\sigma$ is isomorphic to $G$).

\subsection{Solved Problems}

About an eighth of the problems that appeared originally in \cite{MMV} have now been 
solved. We have removed these from the original list. We list them 
here with references for their solutions.
\medskip

\noindent1.11. (J. Martinez \cite{M285}) If ${\cal V}$ is a variety of residually ordered 
$\ell$-groups and satisfies no group laws, then does ${\cal V}$ 
contain a free linearly ordered group? Is ${\cal V}={\cal R}$? 
\medskip

Yes. See N. Ya. Medvedev \cite{M05}.
\medskip

\noindent1.13. (V. M. Kopytov, \cite{K289} Problem 5.23 and \cite{H283} Problem 40)
Is the variety of weakly abelian $\ell$-groups the 
smallest variety of $\ell$-groups that contains all nilpotent $\ell$-groups?  Is the free weakly abelian 
$\ell$-group on $n$ generators reisdually ordered nilpotent?
Residually an ultraproduct of ordered nilpotent groups?
\medskip

The answers to all parts of 1.13 are No. See V. V. Bludov, A. M. W. Glass \cite{BG06}.
\medskip

\noindent1.19. (W. Ch. Holland) If $F(x,y)$ is the free
group and $w(x,y)\in F(x,y)$, is there a fixed-point free image 
$w(\bar x,\bar y)\in Aut(\R)$?
\medskip

Yes. See S. A. Adeleke, W. Ch. Holland \cite{AH94}.
\medskip

\noindent1.24. (A. I. Kokorin \cite{K289} Problem 5.20) Is the elementary theory of lattices of 
$\ell$-ideals of abelian $\ell$-groups decidable?
\medskip

Yes. See N. Ya. Medvedev \cite{M5}.
\medskip

\noindent1.25. (V. M. Kopytov \cite{K289} Problem 5.24) Is any free residually ordered $\ell$-group residually 
ordered soluble?
\medskip

No. See N. Ya. Medvedev \cite{M05}.
\medskip

\noindent1.43. Is it possible to embed a recursively generated $\ell$-group with soluble word problem
in a finitely presented $\ell$-group (the finite presentation being under the group 
and lattice operations)? 
Can a recursively generated $\ell$-group be $\ell$-embedded in a finitely presented $\ell$-group iff it can be defined by a recursively enumerable set of relations?

\medskip

Yes to both questions. Any recursively generated $\ell$-group that is $\ell$-embeddable in a finitely presented $\ell$-group must be defined by a recursively enumerable set of relations. The converse follows from \cite{G06} since every such $\ell$-group can be embedded in a finitely generated $\ell$-group that is defined by a recursively enumerable set of relations (see \cite{G06}). This answers the second question
and hence the first. Also see \cite{G08} for an algebraic characterisation of recursively generated $\ell$-groups with soluble word problem.

\subsection{Comments}

Partial answers have been provided to some of the questions. 
We now list these.
\medskip

\noindent1.5. V. M. Kopytov \cite{K282} proved that the lattice order 
on a free $\ell$-group is extendable to a total order.
\medskip

\noindent1.7. Yes if the abelian group has size at most the continuum, 
since it can then be embedded in $\R$ and hence made into an 
archimedean totally ordered group.
\medskip

\noindent1.10. If ${\cal V}_1={\cal V}_2={\cal A}$ ($\cal A$ the variety of abelian
$\ell$-groups), then ${\cal A}^2$ has a finite basis of laws (A. M. W. Glass, W. Ch. Holland, 
S. H. McCleary \cite{G284}). For other partial cases see \cite{L109}.
\medskip

\noindent1.12. Yes if the Engel class is bounded. See N. Ya. Medvedev \cite{M286}.
\medskip

\noindent1.14. The only {\it varieties} which are idempotents are the trivial 
variety, the variety of all normal-valued $\ell$-groups, and the 
variety of all $\ell$-groups (A. M. W. Glass, W. Ch. Holland, 
S. H. McCleary \cite{G284}).
\medskip

\noindent1.16. The free abelian $\ell$-group on a finite number $n$ of generators has undecidable theory iff $n\geq 3$ (A. M. W. Glass, A. J. Macintyre, 
F. Point \cite{GMP}).
\medskip

\noindent1.17. The free $\ell$-group on $n$ generators has soluble
word problem (W. Ch. Holland, S. H. McCleary \cite{WCH}).
Hence, the subtheory 

\noindent $\{ (\forall x_1,\dots,x_n)(w(x_1,\dots,x_n)=1)\mid w\in F_n\}$
of Th($F_n$) is decidable.
\medskip

\noindent1.18. Finitely generated free abelian $\ell$-groups with distinct ranks
are not elementary equivalent (A. M. W. Glass, A. J. Macintyre, 
F. Point \cite{GMP}).
\medskip

\noindent1.20. For partial results see \cite{C20}. A positive answer implies a positive answer for Problem 1.22.
\medskip

\noindent1.26. The word problem for such $\ell$-groups is soluble (V. M. Kopytov  
\cite{K290}).
\medskip

\noindent1.27. If the variety contains all metabelian $\ell$-groups, then it 
fails the amalgamation property by the example of K. R. Pierce \cite{P972}.
S. A. Gur\-chen\-kov \cite{G103} proved that if a variety of $\ell$-groups $\cal V$ 
satisfies the amalgamation property, then ${\cal V}\subset {\cal R}$. 
\medskip

\noindent1.28. Varieties of $\ell$-groups generated by $Aut(\Omega,\leq)$ were 
described by W. Ch. Holland \cite{H291}.
\medskip

\noindent1.32. The set of all $\ell$-subvarieties of ${\cal R}$ without covers 
(in ${\cal R}$) has at least $6$ elements (N. Ya. Medvedev \cite{M292} and N. Ya. Medvedev, S. V. Morozova \cite{M86}).
\medskip

\noindent1.36. The solubility of the word problem is preserved by free 
products in the variety of all abelian $\ell$-groups and in the variety 
of all $\ell$-groups that are nilpotent class $n$ ($n\in \Z_+$) by V. 
M. Kopytov (see comments on Problem 1.26 above). It is also preserved in the 
variety of all $\ell$-groups (A. M. W. Glass \cite{G006}). 
\medskip

\noindent1.40. Yes for abelian $\ell$-groups (A. Bigard, P. F. Conrad, 
S. Wolfenstein \cite{B294}).
\medskip

\noindent1.41. Yes for abelian finite-valued groups (A. Bigard, P. F. Conrad, 
S. Wolfenstein \cite{B294}).
\medskip

\noindent1.55. Yes if the number of generators is at most $3$ 
(A. M. W. Glass, F. Point \cite{GP}). No if we consider finitely generated 
abelian $\ell$-groups instead of finitely presented abelian $\ell$-groups. 
For let $\xi$ be any irrational real number. One can view the two generator group $D(\xi):=\Z \oplus \Z\xi$ as an o-subgroup of $\R$; so $D(\xi)$ is a dense totally ordered abelian $\ell$-group. But if $G$ and $H$ are dense finitely generated o-subgroups of $\R$, then $G$ and $H$ satisfy the same sentences iff $\dim(G/pG)=\dim(H/pH)$ for all primes $p$ \cite{RK}. 
Hence $D(\xi)$ and $D(\eta)$ are elementarily equivalent for all irrational real numbers $\xi$ and $\eta$. Consequently, $\{ D(\xi)\mid \xi\in \R\setminus \Q\}$ is an uncountable family of two generator elementarily equivalent pairwise non-isomorphic totally ordered abelian $\ell$-groups in the countable language of $\ell$-groups (or totally ordered groups). 
\medskip

\noindent1.56.  Every finitely presented abelian $\ell$-group is projective in the class of abelian $\ell$-groups (W. M. Beynon \cite{Bey}).
\medskip

\noindent1.57. A. Khelif has found an example of a QFA totally ordered group which is not
prime for the theory of totally ordered groups (\cite{AK}, Theorem 7). No such example is known for the theory of groups.

\section{Right-ordered groups}

\subsection{Open Problems}

If $G$ is an $\ell$-group, then the lattice ordering can be extended to a total ordering on $G$ that is preserved by multiplication on the right. Indeed, the lattice ordering is the intersection of all extension right orderings of $G$. Furthermore, every 
right-orderable group can be embedded (as a group) in a lattice-orderable group.
It is therefore natural to extend the study to right-ordered groups.
An important class of right-ordered groups was developed by Paul F. Conrad.
It satisfies the further condition: if $1<f<g$, then $fg^mf^{-1}>g$ for some $m=m(f,g)\in \Z_+$. 
(Indeed, $m$ can always be taken to be $2$ (see \cite{Na}).)
In his honour, such a right ordering is called a {\it Conrad right-ordering}.
A right-orderable group is said to be {\it Conrad right-orderable} if it admits a Conrad right-ordering. 
\medskip

The following notation will be used.
\medskip

If $H$ is a subgroup of a group $G$, then $I(H)$, the {\it isolator subgroup of} $H$, is the smallest subgroup $I\supseteq H$ such that ($n\in \Z_+$, $g\in G$ and $g^n\in I$) implies $g\in I$. 

$G$ is an $R^*$ {\it group} if a (non-empty) product of conjugates of a non-
identity element is never the identity.

A right-orderable group $G$ is an $(OR)^*$ {\it group} if every right 
partial order can be extended to a right order.

A group $G$ is called a {\it unique product group} if for any finite non-empty subsets 
$X,Y\subseteq G$, there exists $g\in G$ that can be written in a uniquely as $g=xy$ ($x\in X$,
$y\in Y$).
\medskip

\noindent2.1$^*$. (V. M. Kopytov, N. Ya. Medvedev \cite{K289} Problem 16.48) Is every $R^*$ 
group right orderable?
\medskip

\noindent2.7. (A. H. Rhemtulla \cite{R297}) A group $G$ is called {\it strictly bounded} if for some 
positive integer $n$, each of its subgroups of the form 
$\langle x^{y^m}\mid m\in\Z\rangle$ ($x,y\in G$) has a generating set of size 
$n$. Does every strictly bounded right-orderable group have a Conrad right 
order?
\medskip

\noindent2.8. (A. H. Rhemtulla \cite{R298}) Does every right-orderable $n$-Engel 
group have a Conrad right order? (cf. Problems 1.12, 2.9, 2.25, 2.29, 3.4, 3.8, 3.9, 3.18)
\medskip

\noindent2.9. Is every torsion-free Engel group right orderable? (cf. Problems 
1.12, 2.8, 2.25, 2.29, 3.4, 3.8, 3.9, 3.18)
\medskip

\noindent2.10$^*$. (A. H. Rhemtulla) For a group $G$, let $F(G)$ be 
the family of all finite non-empty subsets of $G$. Then $F(G)$ is a monoid 
under multiplication induced by the group operation. A homomorphism 
$\sigma: F(G)\rightarrow G$ such that $\{g\}\sigma=g$ for all $g\in G$ is 
called a {\it retraction}. Assume that $G$ has at least one retraction. 
Is $G$ a unique product group?
\medskip

\noindent2.11$^*$. (A. I. Malcev \cite{K289} Problem 1.6) Is the group ring of a right 
orderable group embeddable in a Skew field?
\medskip

\noindent2.12$^*$. Can every metabelian right-ordered (right-orderable, $(OR)^*$) group be o-embedded in a divisible metabelian right-ordered 
(right-orderable, $(OR)^*$) group? The same questions for soluble.
\medskip

\noindent2.13$^*$. Under what conditions is it possible to embed
a right-orderable group in a right-orderable group which has a 
divisible maximal locally nilpotent subgroup?
\medskip

\noindent2.14$^*$. (A. H. Rhemtulla \cite{H283} Problem 20) Is there a finitely generated 
simple right-orderable group? 
\medskip

\noindent2.16. Is there an $(OR)^*$ simple group? 
\medskip

\noindent2.17$^*$. If $G_1,G_2\in (OR)^*$, then is 
$G_1\times G_2\in (OR)^*$? In particular, the special case when $G_2=\Z$?
\medskip

\noindent2.19$^*$. (B. H. Neumann). Is there a universal countable 
right-orderable group? 
\medskip

\noindent2.20. Can the class of right-orderable (locally indicable) 
groups be defined by an independent system of quasi-identities in the class 
of torsion-free groups?
\medskip

\noindent2.23. Describe the right-orderable groups whose set of right 
relatively convex subgroups is a sublattice of the set of all subgroups.
\medskip

\noindent2.24. Describe the right-orderable groups whose set of right 
relatively convex subgroups forms a distributive lattice. 

\begin{center}
{\bf New problems}
\end{center}

\noindent2.25$^*$. Is every right-orderable $n$-Engel group nilpotent?  
(cf. Problems 1.12, 2.8, 2.9, 2.29, 3.4, 3.8, 3.9, 3.18)
\medskip

\noindent2.26. Is there a universal countable locally indicable group?
\medskip

\noindent2.27. (A. H. Rhemtulla) Let $G$ be a (soluble) orderable group.
Suppose that all right orders on $G$ are Conrad right orders.

(i) Is every right order on every subgroup of $G$ a Conrad right order?

(ii) Let $H$ be a subgroup of $G$ such that $G=\mbox{I}(H)$. Is every right order on $H$ a Conrad right order?

(iii) Let $H$ be a subgroup of $G$ with $[G:H]$ finite. Is every right order on $H$ a Conrad right order?

(iv) If $G$ is non-abelian, does it contain a non-commutative free semigroup?

(v) Let $K$ be a central extension of $G$. Is every right order on $K$ a Conrad right order?

(vi) If additionally $G$ is finitely generated, then is $G$ residually torsion-free-nilpotent? 
\medskip

\noindent2.28. (V. V. Bludov, A. H. Rhemtulla) Let $G$ be an orderable 
group such that every order on every subgroup of $G$ is central. Is every right 
order on $G$ a Conrad right order?
\medskip

\noindent2.29. (V. V. Bludov \cite{K289} Problem 16.15) Does the set of Engel elements of a (right) 
orderable group form a subgroup? (cf. Problems 1.12, 2.8, 2.9, 2.25, 3.4, 3.8, 3.9, 3.18)
\medskip

\noindent2.30$^*$. (V. V. Bludov) Let $G$ be a right-orderable group and $N$ be a 
nilpotent subgroup of $G$. Let $P$ be a positive cone for some right order on 
$G$. Let $N^{\#}$ be a minimal divisible completion of $N$ and $P^{\#}$ be the 
positive cone of a right order on $N^{\#}$ which is an extension of the right 
order of $N$ with positive cone $P\cap N$. Let $g\in N^{\#}$ and $\leq^g$ 
denote the right order on $N$ induced by the right order of $N^{\#}$ with positive 
cone $(P^{\#})^g$. Is $\leq^g$ extendable to a right order on $G$?
\medskip

\noindent2.31$^*$. (V. V. Bludov) Let $G$ be a right-orderable group, and $A$ a 
normal abelian subgroup of $G$. Denote by $A^{\#}$ a minimal divisible 
completion of $A$. Let $G\leq G^{\#}$ so that $G^{\#}=A^{\#}G$ and 
$G^{\#}/A^{\#}\cong G/A$. Is $G^{\#}$ right orderable?
\medskip

\noindent2.32. (V .V. Bludov, A. H. Rhemtulla) Let $G$ be a soluble right-orderable group. Can $[G:G']$ be finite?
\medskip

\noindent2.33. (V. V. Bludov, A. M. W. Glass \cite{BG008}) Let $G$ be a right-ordered (right-orderable) group, $a_1,\dots a_k\in G$, and  
$w(x_1,\dots,x_n,a_1,\dots,a_k)=1$ be an equation in $G$.

(i) Is there an  order-embedding of $G$ into a right-ordered group $G^{\ast}$
containing a solution of the equation $w(x_1,\dots,x_n,a_1,\dots,a_k)=1$?

(ii) If $G$ is soluble, is there an order-embedding of $G$ into a right 
ordered soluble group $G^{\ast}$ (of the same solubility length) 
containing a solution of the equation $w(x,a_1,\dots,a_k)=1$?

(iii) If $G$ is a locally indicable group, is there an order-embedding of $G$ into a
right-ordered (locally indicable) group $G^{\ast}$ containing a solution of the
equation $w(x,a_1,\dots,a_k)=1$?
\medskip

\noindent2.34. (P. A. Linnell \cite{K289} Problem 15.49) Is every unique product group right 
orderable?
\medskip

\noindent2.35$^*$. (P. A. Linnell \cite{Lin1}) Can every right-orderable group 
without non-abelian free subgroups be Conrad right-ordered? (cf. Problem~2.27 (iv)) 

\subsection{Solved Problems}

Nearly a third of the original problems in the original list \cite{MMV} have now been solved.
\medskip

\noindent2.2. Is an HNN-extension of a right orderable group right 
orderable? 
\medskip

\noindent2.3. Is an HNN-extension of a Conrad right-orderable group right 
orderable? 
\medskip

No to both questions. See V. V. Bludov, A. M. W. Glass \cite{BG09}.
\medskip

\noindent2.4. (V. V. Bludov) If $G_1,G_2$ are right orderable groups, 
then is $G_1*_H G_2$ right orderable if $H$ is cyclic or abelian rank 1?
\medskip

Yes. See V. V. Bludov, A. M. W. Glass \cite{BG09}.
\medskip

\noindent2.5. (V. V. Bludov \cite{K289}) If $G_1,G_2$ are orderable groups, then 
is $G_1*_H G_2$ right orderable?
\medskip

Yes. See V. V. Bludov, A. M. W. Glass \cite{BG08}.
\medskip

\noindent2.6. If $G_1,G_2$ are right-orderable groups, then 
is $G_1*_H G_2$ right orderable if $H$ is right relatively convex in both factors?
\medskip

Yes. See V. V. Bludov, A. M. W. Glass \cite{BG09}.
\medskip

\noindent2.15. If a non-abelian right-orderable group has no non-trivial normal relatively convex subgroups, then is it simple?
\medskip

No. Consider the right-orderable group 
$G=\langle x,y,z\mid x^2=y^3=z^7=xyz\rangle$ (see G. M. Bergman \cite{Berg1}). 
Since $G$ is finitely generated, it has a maximal proper right relatively convex 
normal subgroup $N$. Since $G=[G,G]$ (see \cite{Berg1}), $G/N$ is non-abelian. Let $\zeta(G)$ be the centre of $G$. It follows from the defining 
relations of $G$, that $x^2,y^3,z^7\in \zeta(G)$. Thus $G=I(\zeta(G))$.
If $\zeta(G)\subseteq N$, then $G$ is a subgroup of $N$ because 
$G=I(\zeta(G))$, a contadiction. Therefore $\zeta(G)\not\subseteq N$.
Hence $G/N$ is non-abelian, has no non-trivial normal right relatively convex 
subgroups, and has a proper non-trival centre.
\medskip

\noindent2.18. (V. M. Kopytov, N. Ya. Medvedev \cite{K289} Problem 16.51) Is there a group 
with exactly a countable infinity of right orders?
\medskip

No. P. A. Linnell \cite{Lin}. (Although there is an error in the proof of Lemma 2.3, this does not affect the proof of Theorem 1.3. which gives a negative answer to 2.18). Subsequently, A. Navas-Flores independently obtained the result from a more dynamical approach when the group is countable (see \cite{Na}), and then in the general case with C. Rivas (see \cite{NRC}).
In contrast, R. N. Buttsworth constructed orderable groups with countably many two-sided orders \cite{B308}.
\medskip

\noindent2.21. (A. H. Rhemtulla) If $G$ is orderable and every 
right order is a Conrad right order, then is $G$ locally nilpotent?
\medskip

No. See V. V. Bludov, A. M. W. Glass and A.H. Rhemtulla \cite{BGR}, Example~4.2, 
Theorem~G, and Proposition~4.3.
\medskip

\noindent2.22. Is there a finitely presented right-orderable 
group with insoluble word problem?
\medskip

Yes. See V. V. Bludov, M. Giraudet, A. M. W. Glass, G. Sabbagh \cite{BGGS}
and V. V. Bludov, A. M. W. Glass \cite{BG09}. It follows now from a result
of V. P. Belkin, V. A. Gorbunov \cite{B318} that the lattice of 
quasi-varieties of right orderable groups contains no maximal elements. 
\medskip

\subsection{Comments}

\noindent2.1. Examples of $R^*$-groups which are not two-sided orderable are
known (\cite{B295}, \cite{R296}, \cite{BL}) but they are all right orderable.
\medskip

\noindent2.10. Some properties of retractable groups were obtained in \cite{B299}.
\medskip

\noindent2.11.  A. I. Malcev \cite{M300} and B. H. Neumann \cite{N301} proved
that the group ring of a two-sided orderable group is embeddable in a Skew field.
\medskip

\noindent2.12. The analogous questions for two-sided orderable and ordered 
groups were solved positively in, respectively, \cite{BM74} and \cite{B003}.
 V. V. Bludov has shown that there are two-sided orderable centre-by-metabelian groups which cannot be 
embedded in any divisible two-sided orderable groups \cite{B005}.
\medskip

\noindent2.13. V. V. Bludov and A. M. W. Glass have shown that this can always be done with abelian in place of 
locally nilpotent  \cite{BG09}.
Moreover, S. A. Gurchenkov has shown that the answer is always positive for two-sided orderable groups  \cite{G259}.
\medskip

\noindent2.14. Examples of perfect ($G=[G,G]$) finitely generated right-orderable groups
were found by G. M. Bergman \cite{Berg1}.
\medskip

\noindent2.17. The analogous question for $O^*$-groups was solved positively 
by A. I. Kokorin \cite{K306} and M. I. Kargapolov \cite{K307} (see also \cite{KK72}).
\medskip

\noindent2.19. There is a universal finitely presented right-orderable group, and ditto for right ordered (V. V. Bludov, A. M. W. Glass \cite{BG09}).
If we ask the same question for finitely generated right-order\underline{ed} groups or finitely generated two-sided totally order\underline{ed} abelian groups, the answer is trivially No: Let $\xi$ be any irrational real number and let $D(\xi)$ be the totally ordered abelian $\ell$-group described in the comments to Problem 1.54. 
If there were a universal countable right-ordered group, say $U$, then we could o-embed $D(\xi)$ in $U$. Let $u(\xi),v(\xi)$ be the images of $1$ and $\xi$ in $U$.
As $\xi$ varies, since $U$ is countable, there is an uncountabe set $X$ such that
$u(\xi)=u(\eta)$ for all $\xi,\eta\in X$.
But for distinct $\xi,\eta\in X$, $\langle u(\xi),v(\xi)\rangle$ and $\langle u(\eta),v(\eta)\rangle$ are not isomorphic, so $v(\xi)\neq v(\eta)$. Thus $U$ contains continuum many elements $\{ v(\xi)\mid \xi\in X\}$, a contradiction. The same argument applies to show that there is no universal countable totally ordered group.
\medskip

We do not know if there is a universal countable totally order\underline{able} group.
\medskip

\noindent2.21. The converse is true (A. H. Rhemtulla \cite{R72} and 
J. S. Ault \cite{A72}). In an e-mail to us, A. Navas has asked the weaker question: If $G$ is right orderable and every 
right order on $G$ is a Conrad right order, then is $G$ locally residually nilpotent?
\medskip

\noindent2.25. Yes if $n\leq 4$. See P. Longobardi, M. Maj \cite{LM98}
and G. Havas, M. R. Vaughan-Lee \cite{HV05} for an arbitrary $4$-Engel group.  
\medskip

\noindent2.30.  A positive solution to this problem implies a positive solution 
of Problem~2.13 (see V. V. Bludov, A. M. W. Glass \cite{BG09}).
\medskip

\noindent2.31. This is true for orderable groups (C. D. Fox \cite{F303}). Although
one can embed a right-orderable group in a right orderable group
with divisible abelian subgroup $A^{\#}$ (see comments to Problem~2.13),  the
construction in \cite{BG09} does not respect the normality of $A^{\#}$.
\medskip

\noindent 2.35.
Every right-orderable amenable group is Conrad right-orderable (D. Witte Morris  \cite{M-W}).
A. Navas has pointed out that there seems to be no known example of a 
right-orderable group satisfying a group identity which is not Conrad right-orderable.
\medskip

\section{Totally ordered groups}

Finally, we consider questions involving the class of ordered (or orderable) 
groups. These are special cases of both $\ell$-groups (with total order) and 
right-ordered groups (where the order is also preserved by multiplication on the left).
\medskip 

Further notation: A group $G$ is said to be an $O^*$-group if every 
partial order on $G$ can be extended to a total order on $G$ making 
$G$ a totally ordered group (or o-group, for short).

A subgroup $H$ of a group $G$ is called {\it infrainvariant} if for each 
$g\in G$, either $H\subseteq g^{-1}Hg$ or $g^{-1}Hg\subseteq H$.
\medskip

\noindent3.1$^*$. (V. V. Bludov) If $G\in R^*$ has a single defining relation, 
then is $G$ orderable?
\medskip

\noindent3.2. Are Bergman's conditions for a free product with 
amalgamation to be orderable also sufficient? (see G. M. Bergman \cite{B208}).
\medskip

\noindent3.3$^*$. (V. M. Kopytov, N. Ya. Medvedev \cite{K289} Problem 16.49) Is the free 
product of $R^*$ groups also $R^*$ ?
\medskip

\noindent3.4$^*$. (A. I. Kokorin \cite{K289} Problem 2.24) Is every torsion-free $n$-Engel group 
orderable? (cf. Problems 1.12, 2.8, 2.9, 2.25, 3.4, 3.8, 3.9, 3.18)
\medskip

\noindent3.6. If $F$ is a free group and $H\lhd F$ with $F/H$ right 
orderable, then is $F/[H,H]$ orderable?
\medskip

\noindent3.7. (A. H. Rhemtulla \cite{R260}) If $F$ is a free group and $H\lhd F$ with $F/H$ a 
periodic $2$-group, then can every $F$-invariant order on $H$ be extended to 
make $F$ an o-group?
\medskip

\noindent3.8. (A. H. Rhemtulla) If $F$ is a free group and $H\lhd F$ with $F/H$ a 
periodic Engel group, then can every $F$-invariant order on $H$ be extended to 
make $F$ an o-group? (cf. Problems 1.12, 2.8, 2.9, 2.25, 3.4, 3.8, 3.9, 3.18)
\medskip

\noindent3.9. (A. H. Rhemtulla \cite{R260}) If $G$ is a torsion-free group and $H$ is a soluble normal 
subgroup of $G$ with $G/H$ a periodic Engel group, then can every $G$-invariant order on $H$ be 
extended to make $G$ an o-group? (cf. Problems 1.12, 2.8, 2.9, 2.25, 3.4, 3.8, 3.9, 3.18)
\medskip

\noindent3.10. If an orderable group has finite special rank, 
then is it soluble? 
\medskip

\noindent3.11. (A. H. Rhemtulla \cite{R297}) If an orderable group has finite 
special rank, are all its convex subgroups isolators of finitely 
generated subgroups?
\medskip

\noindent3.15$^*$. (A. I. Kokorin \cite{K289} Problem 2.28) Can every orderable 
group be embedded in an $O^*$-group?
\medskip

\noindent3.16$^*$. (L. Fuchs, A. I. Malcev \cite{K289} Problem 1.35 written in 
by M. I. Kargapolov) Is there a simple $O^*$-group?
\medskip

\noindent3.17$^*$. (A.I. Kokorin \cite{K289} Problem 3.20) Is the class of 
orderable groups the axiomatic closure of the class $O^*$?
\medskip

\noindent3.18. Are Engel orderable groups $O^*$? (cf. Problems 1.12, 2.8, 2.9, 2.25, 3.4, 3.8, 3.9, 3.18)
\medskip

\noindent3.19. (A.I. Kokorin \cite{KK72} Problem 19) Is every orderable group $H$ embeddable in an orderable 
group $G$ so that any $G$-invariant partial order on $H$ has an extention to 
a total order on $G$?
\medskip

\noindent3.20$^*$. Is every orderable abelian-by-nilpotent group an $O^*$-group?
\medskip

\noindent3.21. Let ${\cal W}^*$ be the class of orderable groups $G$ 
such that every $G$-invariant weakly abelian partial order can be extended 
to a weakly abelian total order. 
Is ${\cal W}^*\setminus O^*\neq\emptyset$?
\medskip

\noindent3.22$^*$. Is Chehata's group an $O^*$-group?
\medskip

\noindent3.23$^*$. Describe $Aut(Dlab[0,1])$ and $Aut(Dlab(\R))$.
\medskip

\noindent3.24$^*$. For each $k$, construct a group $G_k$ that has 
exactly $2k$ orders.  (cf. Problem 1.6)
\medskip

\noindent3.25$^*$. (S. H. McCleary \cite{K289} Problem 12.51) Let $F$ be the free group of finite 
rank $\geq 2$. Does there exist a finite subset $S\subset F$ such that there is a 
unique total order on $F$ under which all elements from $S$ are positive?
\medskip

\noindent3.26. (A. I. Kokorin \cite{KK72} Problem 28) If the centraliser of a relatively 
convex subgroup is $R^*$-isolated, then is the centraliser relatively convex?
\medskip

\noindent3.27. (A. I. Kokorin \cite{KK72} Problem 6) Let $\{ G_i\mid i\in I\}$ be a family of relatively 
convex subgroups of an orderable group $G$, and 
$H=\bigcup_{i\in I}\,G_i$. If $H$ is an infrainvariant subgroup of 
$G$, is it relatively convex?
\medskip

\noindent3.28$^*$. (A. I. Kokorin \cite{K289} Problem 1.46) Find conditions for the 
normaliser of a relatively convex subgroup to be relatively convex.
\medskip

\noindent3.29$^*$. Is there a finitely presented orderable group with insoluble 
word problem? (cf. Problem 2.22).

\begin{center}
{\bf New problems}
\end{center}

\noindent3.30. (V. V. Bludov \cite{B74}) Does there exist an $O^*$-group with
only two orderings: an ordering $\leq$ and its reverse?
\medskip

\noindent3.31$^*$. (A. M. W. Glass \cite{H283} Problem 1$^\prime$) Does there exist a totally ordered 
group all of whose strictly positive elements are conjugate? 
\medskip

\noindent3.32$^*$. (V. V. Bludov) Let $G$ be an abelian-by-nilpotent ordered group. 

(i) Is $G$ order-embeddable in a divisible ordered group?

(ii) If $G$ is also nilpotent-by-abelian, is $G$ order-embeddable in a 
divisible ordered group?
\medskip

\noindent3.33. (V. V. Bludov) Let $G$ be an ordered group $G$ that can be order-embedded in a 
divisible ordered group. We say that $G^{\ast}$ is a {\it minimal divisible completion 
of} $G$ if $G^{\ast}=\mbox{I}(G)$. Let ${\cal D}_G$ be the class of all 
minimal divisible completions of $G$. 

(i) Does there exist a group $K\in{\cal D}_G$ such that for any 
$H\in{\cal D}_G$ there is a homomorphism of $K$ onto $H$ that is the identity on $G$?

(ii) The same question for metabelian $G$.

\subsection{Solved Problems}

Again, about an eighth of the problems in the original list have now been 
solved.
\medskip

\noindent3.5. Does there exist a soluble orderable group $G$ with 
$G/[G,G]$ periodic?
\medskip

Yes (V. V. Bludov, V. M. Kopytov, A. H. Rhemtulla \cite{BKR}).
\medskip

\noindent3.12. If an orderable group contains no non-abelian free 
subsemigroups, then is it nilpotent?
\medskip

\noindent3.13. If a residually nilpotent orderable group contains no 
non-abelian free subsemigroups, then is it nilpotent?
\medskip

No to both questions as can be seen by Golod examples (V. V. Bludov, 
A. M. W. Glass, A. H. Rhemtulla \cite{BGR}).
\medskip

\noindent3.14. If $G$ is an orderable nilpotent-by-abelian group, 
can it be embedded in a divisible such? The same question with ordered 
in place of orderable and the embedding an o-embedding.
\medskip

No to both parts (V. V. Bludov \cite{B005}).
\medskip

\subsection{Comments}

\noindent3.1. A torsion-free group with a single defining relation admits a
Conrad right order (S. D. Brodsky \cite{B194}).
\medskip

\noindent3.3. Yes if one of the factors is right orderable (N. Ya. Medvedev 
\cite{M264}).
\medskip

\noindent3.4. Yes if it is lattice-orderable (see Problem 1.12).
\medskip

\noindent3.15. Free groups are embeddable in $O^*$-groups (V. M. Kopytov 
\cite{K179}).
\medskip

\noindent3.16. There are simple groups in which all partial orders are extendable to
lattice orders (N. Ya. Medvedev \cite{M160}).
\medskip

\noindent3.17. The class of right-orderable groups does not coincide with the
axiomatic closure of the class $(OR)^*$ (N.Ya. Medvedev \cite{M214}).
\medskip

\noindent3.20. Orderable metabelian groups are $O^*$-groups (A. I. Kokorin
\cite{K313}) and abelian-by-nilpotent $R^*$-groups which are nilpotent-by-abelian are also $O^*$-groups (P. Longobardi, M. Maj, A. H. Rhemtulla \cite{LMR}). Orderable groups with nilpotent commutator subgroup need not be $O^*$-groups
(M. I. Kargapolov \cite{K307}).
\medskip

\noindent3.22 and 3.23. Definitions and basic properties of Chehata's and
Dlab's groups can be found in \cite{K4} and \cite{K7}.
\medskip

\noindent3.24. If $k=4n$, then there is a metabelian group $G_k$ with exactly 
$2k$ orders (V. M. Kopytov \cite{K314}). In \cite{M253}, N. Ya. Medvedev 
constructed groups $G_k$ for $k=2(4n+3)$. But 
the existence of a group $G_{10}$ (with exactly $10$ total orders) and 
$G_{18}$ with exactly $18$, for example, remain open.
\medskip

\noindent3.25. There is no non-empty finite subset $S\subset F$ with a 
unique right order on $F$ under which all elements from $S$ are strictly positive (S. H. McCleary \cite{M316}).
\medskip

\noindent3.28. The normaliser of a relatively convex subgroup
need not be relatively convex (M. I. Kargapolov, A. I. Kokorin, V. M. Kopytov
\cite{K317}).
\medskip

\noindent3.29. If there is a finitely presented orderable group with insoluble 
word problem, then the lattice of quasi-varieties of orderable groups 
contains no maximal elements (see V. P. Belkin, V. A. Gorbunov \cite{B318}).
\medskip

\noindent3.31. A positive solution to this problem would imply positive solutions
for Problems 3.16 and 3.30.
\medskip

\noindent3.32. V. V. Bludov has shown that the answer is yes if $G$ is metabelian \cite{B003}
and no for arbitrary ordered groups \cite{B005}. Also see the comments to Problem 3.14.
\bigskip

{\bf Acknowledgements}
\medskip

This research was supported by grants from the Royal Society, the London Mathematical Society (Scheme IV, Travel) and Queens' College, Cambridge. 
We are most grateful to them for facilitating this work and to the College and DPMMS for their hospitality.

Authors' addresses:
\bigskip

\noindent\parbox[t]{6cm}{
V. V. Bludov:

\medskip
 
Department of Mathematics, Physics, 
and Informatics,
 
Irkutsk State Teachers 
Training University,

Irkutsk  664011,

Russia

\medskip

vasily-bludov@yandex.ru

\bigskip

V. M. Kopytov:
\medskip

Institute of Mathematics,

Academician Koptyug Avenue, 4,

Novosibirsk, 630090

Russia.
\medskip

kopytov@academ.org
}\hfill\parbox[t]{6.1cm}{
A. M. W. Glass: 

\medskip

Queens' College,

Cambridge CB3 9ET,

England
\medskip

and 
\medskip

Department of Pure Mathematics and 
Mathematical Statistics,
 
Centre for Mathematical Sciences, 

Wilberforce Rd., 

Cambridge CB3 0WB, 

England

\medskip

amwg@dpmms.cam.ac.uk
}

\end{document}